\documentclass[9pt,amsfonts]{article}
\usepackage{epsf,pgf,etex,graphicx,amssymb,eucal,mathrsfs}
\usepackage{latexsym,amsmath,epsfig,epic,eepic,xcolor}
\usepackage[all]{xy}
\usepackage{wasysym,ntheorem}
\usepackage{epsf,graphicx,pgf,etex,amscd,pstricks}

\usepackage[hypertexnames=false,backref=page,pdftex,
 	pdfpagemode=UseNone,
 	breaklinks=true,
 	extension=pdf,
 	colorlinks=true,
 	linkcolor=blue,
 	citecolor=red,
 	urlcolor=blue,
 ]{hyperref}

\newcommand{\R}{\mathbb{R}}

\newcommand{\C}{\mathbb{C}}

\newcommand{\HR}{\mathscr{H}}

\newcommand{\HHC}{\mathbb{H}^{\infty}_{\mathbb{C}}}
\newcommand{\ID}{{\rm Id}}
\newcommand{\SUD}{{\rm SU}(2)}
\newcommand{\PON}{{\rm PO}(1,n)}

\newtheorem{theorem}{{\bf Theorem}}

\newtheorem{prop}[theorem]{Proposition}

\theoremnumbering{Alph}
\newtheorem*{main}{{\bf  Main Theorem}}

\title{Hyperbolic spaces, principal series and ${\rm O}(2,\infty)$}
\author{Pierre Py and Arturo S\'anchez\footnote{The authors were partially supported by the french project ANR AGIRA, by project papiit IA100917 from DGAPA UNAM and by the CNRS UMI 2001 (Laboratoire Solomon Lefschetz LaSol, Cuernavaca, M\'exico). Arturo S\'anchez was also supported by a CONACYT fellowship (no.~701896).}}
\date{November 2018}

\begin{document}

\maketitle

\begin{abstract} 
We prove that there exists no irreducible representation of the identity component of the isometry group ${\rm PO}(1,n)$ of the real hyperbolic space of dimension $n$ into the group ${\rm O}(2,\infty)$, if $n\ge 3$. This is motivated by the existence of irreducible representations (arising from the spherical principal series) of ${\rm PO}(1,n)^{\circ}$ into the groups ${\rm O}(p,\infty)$ for other values of $p$. 
\end{abstract}


\section{Introduction}

\subsection{Strongly nondegenerate bilinear forms of finite index}
Consider a separable Hilbert space $\mathscr{H}$, and a Hilbert basis $(e_{i})_{i\ge 1}$ of $\mathscr{H}$. Define a bilinear form $B_{p}$ on $\mathscr{H}$ by the formula
$$B_{p}(x,x)=\sum_{j=1}^{p}x_{i}^{2}-\sum_{j\ge p+1}x_{j}^{2},$$
where $x$ is an arbitrary vector in $\mathscr{H}$ written as $x=\sum_{j=1}^{+\infty}x_{j}e_{j}$. The forms $\pm B_{p}$ can be characterized intrinsically as the unique (up to isomorphism) strongly nondegenerate bilinear forms of index $p$ on a separable Hilbert space, see~\cite{bim}. The word {\it strongly} here refers to a completeness condition \cite[\S 2]{bim}. The isometry group of $B_{p}$ is denoted by ${\rm O}(p,\infty)$. It is interesting to study which locally compact groups admit irreducible representations into this group\footnote{All the representations that we consider in this text will be continuous, i.e. orbit maps $g\mapsto g\cdot v$ are continuous.}. This question was suggested by Gromov~\cite[\S 6]{gromov} and studied by Duchesne who constructed tools to establish superigidity results with target ${\rm O}(p,\infty)$~\cite{duchesne0,duchesne1,duchesne15}. See also~\cite{dlp} for more recent results in this direction. Note also that the groups ${\rm O}(p,\infty)$ were already studied in the 40's and 60's by Ismagilov, Naimark and Pontryagin, see e.g.~\cite{ism,ism2} as well as the references cited in~\cite{monodpy}. 

A natural family of groups admitting interesting irreducible representations into the groups ${\rm O}(p,\infty)$ is the family of isometry groups of the real hyperbolic spaces. This fact has been well-known to representations theorists for years~\cite{johnsonwallach,sally}, and was put in a geometric context more recently~\cite{delzantpy,monodpy}. We describe this in the next section. 

\subsection{The spherical principal series and representations into ${\rm O}(2,\infty)$}

Fix $n\ge 2$ and let $\mathbb{H}^{n}$ be the $n$-dimensional real hyperbolic space. We identify its isometry group with the group ${\rm PO}(1,n)={\rm O}(1,n)/\{\pm \ID\}$. The irreducible representations of $\PON$ into the groups ${\rm O}(p,\infty)$ that we alluded to above come for the study of the so-called {\it spherical principal series}. We recall now the definition of this classical representation theoretic object.

We fix a point $o\in \mathbb{H}^{n}$ and denote by $K\subset \PON$ its stabilizer and by $\mu$ the unique $K$-invariant volume form of total volume $1$ on the boundary $\partial \mathbb{H}^{n}$ of $\mathbb{H}^{n}$. The spherical principal series is a family $(\pi_{s})_{s\in \mathbb{C}}$ of representations of ${\rm PO}(1,n)$ on the space ${\rm L}^2(\partial \mathbb{H}^{n},\mu)$. It is defined by the formula:
$$\pi_{s}(g)(f)=\vert {\rm Jac}(g^{-1})\vert^{\frac{1}{2}+s}\cdot f\circ g^{-1}\;\;\;\;\; (g\in {\rm PO}(1,n), f\in {\rm L}^2(\partial \mathbb{H}^{n},\mu)).$$
Here ${\rm Jac} (g)$ is the Radon-Nikodym derivative of the measure $\mu$ with respect to an element $g\in {\rm PO}(1,n)$.
The study of the properties of $\pi_{s}$ depending on the value of $s$ is a classical topic. We will not describe all the known properties of these representations but we only mention a few classical facts. If $s$ is purely imaginary, $\pi_{s}$ is unitary. Also, the representation $\pi_{s}$ is irreducible if and only if $s\notin\{\pm(\frac{1}{2}+\frac{k}{n-1}), k\in \mathbb{N}\}$~\cite{wallach}. When $s$ is real and positive, it is known that there exists a continuous invariant bilinear form $B_{s}$ for $\pi_{s}$. This form is nondegenerate and of finite index when $\pi_{s}$ is irreducible. We call $p(s)$ its index, which is locally constant on $\R_{+}^{\ast}-\{\frac{1}{2}+\frac{k}{n-1}, k\in \mathbb{N}\}$. Note that the bilinear form $B_{s}$ is not strongly non-degenerate, but it can be naturally completed, thus providing an irreducible representation that we shall call $\overline{\pi}_{s}$:
$$\overline{\pi}_{s} : \PON \to {\rm O}(p(s),\infty).$$

\noindent The nonzero values of the function $s \mapsto p(s)$ (when $s>0$ and $s\notin \{\frac{1}{2}+\frac{k}{n-1}\}$) can be checked to be the set of integers of the form
\begin{equation}\label{indicespossibles}
\binom{n-1+j}{n-1}, j\ge 0.
\end{equation}
\noindent 
We refer to~\cite[\S 3.B]{monodpy} for a detailed discussion of this fact. This leads naturally to the following question (where we restrict to the identity component of the isometry group of $\mathbb{H}^{n}$, for simplicity).

 {\it Is it true that any irreducible representation $\PON^{\circ} \to {\rm O}(p,\infty)$ is conjugated to one of the representations $\overline{\pi}_{s}$?} 

\noindent A more specific question goes as follows. 

{\it If $p$ is not equal to one of the integers of the form~\eqref{indicespossibles}, is it true that there is no irreducible representation of $\PON^{\circ}$ into ${\rm O}(p,\infty)$?} 

\noindent Recall that the first gap occuring in the values of the indices is the set $\{2, \ldots , n-1\}$. This problem was studied in~\cite{monodpy} where the authors proved the following:

\begin{theorem} Let $p$ be an integer with $2<p<n$ where $n>4$. Then there is no irreducible representation ${\rm PO}(1,n)^{\circ}\to {\rm O}(p,\infty)$.
\end{theorem}

Assuming that $p\in \{2,\ldots , n-1\}$, this left open the problem of the existence of irreducible representations ${\rm PO}(1,n)^{\circ}\to {\rm O}(p,\infty)$ in the following cases: $p=2$, $n\ge 3$ and $p=3$, $n=4$. The main result of this note is the following theorem, which settles the case $p=2$. 

\begin{main}\label{mainth} Let $n\ge 3$. Then there is no irreducible representation 
\begin{equation}\label{truc} 
{\rm PO}(n,1)^{\circ}\to {\rm O}(2,\infty).
\end{equation}
\end{main}

Let us describe the strategy of the proof. Using a result from~\cite{monodpy}, we prove that if $\varrho$ is an irreducible representation as in~\eqref{truc}, then its complexification cannot be irreducible. One is thus led to study closed invariant complex subspaces in the complexification. We observe that the complexified Hilbert space, possibly modulo a negative definite invariant factor, is always the direct sum of at most two closed invariant irreducible strongly nondegenerate subspaces. Looking at the signature of these subspaces we are thus led to consider representations into the groups ${\rm U}(2,\infty)$ and ${\rm U}(1,\infty)$, which appear as subrepresentations of the complexification of $\varrho$. From these considerations, we finally derive a contradiction.

In section~\ref{comp}, we make a few general observations about linear representations of arbitrary groups  into ${\rm O}(2,\infty)$, assuming that their complexification is not irreducible, before proving our main theorem in section~\ref{proof}. It is possible that some of our ideas could be used to disprove the existence of irreducible representations into ${\rm O}(2,\infty)$ for other groups. 

{\bf Acknowledgements.} We would like to thank Bruno Duchesne and Nicolas Monod for their comments on this note. 


\section{Complexifications of representations in ${\rm O}(2,\infty)$}\label{comp}

In this section we consider a topological group $G$, a Hilbert space $\mathscr{H}$ and a strongly nondegenerate bilinear form $B$ on $\HR$, of signature $(2,\infty)$. We denote by ${\rm O}(\mathscr{H},B)$ the group of all bijective linear operators of $\mathscr{H}$ preserving $B$ (such an operator is automatically bounded). Let
$$\varrho : G \to {\rm O}(\mathscr{H},B)$$
be an irreducible representation. We let $\varrho_{\C}$ be the complexification of $\varrho$. This is a linear representation of $G$ on the space $\HR \otimes \C$. We identify $\HR$ with $\HR\otimes 1\subset \HR\otimes \C$ and denote by $Re, Im$ the maps $\HR\otimes \C\to \HR$ which send $x+iy$ ($x,y\in \HR$) to $x$ and $y$ respectively. We still denote by $B$ the Hermitian extension of $B$ to $\HR\otimes \C$. 

We assume that $\varrho_{\C}$ is {\it not} irreducible. This means that there exists a closed, non-trivial $G$-invariant complex subspace in $\HR\otimes \C$. In the following proposition we discuss properties of any such complex subspace. This proposition remains true whenever $B$ is strongly nondegenerate of any finite index (not necessarily equal to $2$). 

\begin{prop}\label{blablabla} Let $W\subset \HR\otimes \C$ be a nonzero, proper closed $G$-invariant complex subspace. 
\begin{enumerate}
\item The restriction of $Re$ and $Im$ to $W$ are injective. 
\item The dimension of $W$ is infinite. 
\item The restriction of $B$ to $W$ is strongly non-degenerate. 
\end{enumerate}
\end{prop}
{\it Proof.}  If the restriction of $Im$ to $W$ is not injective, $W\cap \HR$ is a nonzero closed $G$-invariant subspace of $\HR$, hence must be equal to all of $\HR$ since $\varrho$ is irreducible. Hence $\HR\subset W$ and since $W$ is complex $W=\HR\otimes \C$, which is a contradiction. The proof is similar for $Re$. 

 If the dimension of $W$ was finite, $Re(W)$ and $Im(W)$ would be finite dimensional $G$-invariant subspaces of $\HR$, hence they should both be equal to $\{0\}$. This would imply $W=\{0\}$, a contradiction. Hence $W$ is infinite dimensional. 

 Let $N_{W}$ be the radical of the restriction of $B$ to $W$. This space is isotropic, hence finite dimensional since $B$ has finite index. Since it is also $G$-invariant, it must be $\{0\}$ by the previous item. Hence the restriction of $B$ to $W$ is non-degenerate. The fact that this restriction is strongly nondegenerate now follows from Proposition 2.8 in~\cite{bim}. Note that all results in~\cite{bim} are stated in the real case but also hold in the complex case.\hfill $\Box$

The next proposition is a particular case of a result due to Ismagilov~\cite{ism} (see also~\cite[p. 154]{sas} for a detailed proof). The proof essentially consists in observing that, thanks to Zorn's lemma, there exists a {\it maximal} invariant closed complex subspace of  $\HR \otimes \C$ on which $B$ is negative definite. Combining this observation with the fact that $B$ has finite index and with Proposition~\ref{blablabla}, one obtains the result. 

\begin{prop}\label{dichotomy} There is a $G$-invariant orthogonal decomposition into closed strongly nondegenerate complex subspaces of one of the following type.
\begin{enumerate}
\item $\HR \otimes \C=V\oplus W$ where $W$ is negative definite, $V$ is of signature $(2,\infty)$ and the action of $G$ on $V$ is irreducible.
\item $\HR \otimes \C=V_{1}\oplus V_{2}\oplus W$ where $W$ is negative definite, each $V_{i}$ has signature $(1,\infty)$ and the action of $G$ on each $V_{i}$ is irreducible. 
\end{enumerate} \end{prop}

In the next section we turn to the special case where $G={\rm PO}(1,n)^{\circ}$.


\section{Representations of ${\rm PO}(1,n)^{\circ}$ into ${\rm O}(2,\infty)$}\label{proof}

We now turn to the proof of our main theorem. We consider an irreducible representation 
$$\varrho : {\rm PO}^{\circ}(1,n)\to {\rm O}(2,\infty)$$
with $n\ge 3$ and we are looking for a contradiction. As before, we call $\mathscr{H}$ the Hilbert space and $B$ the bilinear form underlying the representation $\varrho$. We fix once and for all a maximal compact subgroup $$K\subset {\rm PO}^{\circ}(1,n).$$
Note that $K$ is isomorphic to the group ${\rm SO}(n)$. We start with the following elementary observation, already used in~\cite{monodpy}. 

\begin{prop}\label{pointwf} There exists a $2$-dimensional, positive definite subspace $P\subset \mathscr{H}$ which is pointwise fixed by $K$. 
\end{prop}

\noindent {\it Proof.} Since this simple result will be used several times here, we recall its proof, although it already appears  in the proof of Theorem 5.3 in~\cite{monodpy}. Since $K$ is compact, it must fix a point in the symmetric space of ${\rm O}(2,\infty)$ which can be identified with the space of positive definite planes in $\mathscr{H}$. Hence there exists a positive definite $K$-invariant plane $P\subset \mathscr{H}$. The action of $K$ on $P$ is given, up to conjugacy, by a homomorphism $K\to {\rm O}(2)$. Any such morphism is trivial hence $P\subset \mathscr{H}^{K}$.\hfill $\Box$

The following statement is proved in~\cite[Prop. 5.4]{monodpy}. Since we will apply it repeatedly to the pair $({\rm PO}^{\circ}(1,n),{\rm SO}(n))$, we state it explicitly here. 

\begin{prop}\label{gelfand} Let $(G,K)$ be a Gelfand pair. Let $\pi$ be a linear representation of $G$ on a complex Hilbert space $W$ preserving a continuous, strongly nondegenerate sesquilinear form of finite index. If $\pi$ is irreducible, then the space $W^{K}$ of $K$-invariant vectors has complex dimension at most $1$. 
\end{prop}

The previous two propositions imply that the complexification $\varrho_{\C}$ of $\varrho$ cannot be irreducible. Indeed the real dimension of the space of $K$-fixed points in $\HR$ is equal to the complex dimension of the space of $K$-fixed points in $\HR\otimes \C$:
$${\rm dim}_{\mathbb{R}}\mathscr{H}^{K}={\rm dim}_{\mathbb{C}} (\mathscr{H}\otimes \mathbb{C})^{K}.$$
If $\varrho_{\C}$ is irreducible, the right-hand side is less or equal than $1$ according to Proposition~\ref{gelfand}, whereas the left-hand side is greater or equal to $2$ according to Proposition~\ref{pointwf}. Hence $\varrho_{\C}$ cannot be irreducible. We will now apply Proposition~\ref{dichotomy}. In the next two subsections, we deal separately with the two cases appearing in that proposition.

\subsection{Irreducible representations in ${\rm U}(2,\infty)$}

We assume here that we have a ${\rm PO}(1,n)^{\circ}$-invariant orthogonal decomposition into closed strongly nondegenerate complex subspaces
$$\HR\otimes \C = V\oplus W$$
where $V$ is irreducible of signature $(2,\infty)$ and where $W$ is negative definite. Consider the space of positive definite complex $2$-dimensional planes in $V$. One shows exactly as in Proposition~\ref{pointwf} that the group $K$ must fix a point $P$ in that space. Consider now the action of $K$ on $P$. Since this action preserves the restriction of $B$ to $P$, this defines a homomorphism from $K$ to the unitary group ${\rm U}(2)$ (well-defined up to conjugacy). But one has the following proposition. 

\begin{prop}\label{matrices} Let $n\ge 3$. Any continuous homomorphism from ${\rm SO}(n)$ to ${\rm U}(2)$ is trivial. 
\end{prop}

This yields a contradiction: indeed the space $V^{K}$ must contain $P$ according to the above proposition. Hence ${\rm dim}_{\C}V^{K}$ is greater or equal to $2$, which contradicts again Proposition~\ref{gelfand} applied this time to the complex representation of ${\rm PO}(n,1)^{\circ}$ on $V$. 

We now prove Proposition~\ref{matrices}. 

\noindent {\it Proof of Proposition~\ref{matrices}.} Let $\varphi : {\rm SO}(n)\to {\rm U}(2)$ be a continuous (hence smooth) homomorphism. The image of $\varphi$ must be contained in ${\rm SU}(2)$ so we think of $\varphi$ as a map from ${\rm SO}(n)$ to ${\rm SU}(2)$. For $n=3$, the Lie algebra of ${\rm SO}(3)$ being simple, $\varphi$ must be a local isomorphism if it is not trivial. But this implies that $\varphi$ is a covering map. The group ${\rm SU}(2)$ being simply connected,  $\varphi$ must be an isomorphism. This is a contradiction since ${\rm SO}(3)$ is not simply connected. (Of course, it is classical that there is a local isomorphism going the other way around, i.e. a $2$-sheeted cover $\SUD \to {\rm SO}(3)$.)

We now deal with the case $n=4$. Recall that there is a $2$-sheeted cover $\pi : \SUD \times \SUD \to {\rm SO}(4)$ whose kernel is generated by $(-{\rm Id},-{\rm Id})$. We consider the composition 
$$\varphi \circ \pi : \SUD \times \SUD \to \SUD.$$
Its restriction to each factor $\SUD \times \{1\}$ and $\{1\}\times \SUD$ is either trivial or an isomorphism onto $\SUD$. Since the two factors commute, it cannot be an isomorphism on each factor since $\SUD$ is not abelian. Hence the map $\varphi \circ \pi$ factors through one of the projections from $\SUD\times \SUD$ onto one of its factor; for instance the first one. Hence we can write:
\begin{equation}\label{abc}
\varphi \circ \pi= \Phi \circ p_{1}
\end{equation}
where $p_1$ is the first projection $\SUD\times \SUD\to \SUD$ and $\Phi : \SUD\to \SUD$ is an isomorphism (if $\varphi$ is nontrivial). Applying this identity to $(-{\rm Id},-{\rm Id})$ one sees that $-{\rm Id}$ lies in the kernel of $\Phi$. This contradicts the fact that $\Phi$ is an isomorphism. Hence $\varphi$ is trivial. 

Finally, for $n\ge 5$, the Lie algebra of ${\rm SO}(n)$ being simple of dimension $>3$, the morphism $\varphi$ must also be trivial.\hfill $\Box$

\subsection{Actions on complex hyperbolic spaces}

To study complex representations of ${\rm PO}(1,n)^{\circ}$ into the group ${\rm U}(1,\infty)$, we will need the following proposition. 

\begin{prop}\label{totallyreal} Let $\alpha : {\rm PO}(1,n)^{\circ}\to {\rm U}(1,\infty)$ be an irreducible representation. Then $\alpha$ is the complexification of an irreducible representation ${\rm PO}(1,n)^{\circ} \to {\rm O}(1,\infty)$.  
\end{prop}

In~\cite{burgeriozzi} Burger and Iozzi classified representations of finitely generated groups into the group ${\rm PU}(1,m)$ with vanishing bounded K\"ahler class.  Here we need to understand complex representations of the group ${\rm PO}(1,n)^{\circ}$ into the group ${\rm U}(1,\infty)$. They are automatically of zero bounded K\"ahler class since the second continuous bounded cohomology group of ${\rm PO}(1,n)^{\circ}$ is zero. So one could try to establish a generalization of Burger and Iozzi's work~\cite{burgeriozzi}where the ``source group" is non-discrete and the target is the isometry group of the infinite dimensonial complex hyperbolic space. We will not pursue this objective here but will only study the case where the source group is ${\rm PO}(1,n)^{\circ}$, which is much simpler and sufficient to establish Proposition~\ref{totallyreal}.   

We first explain how to conclude the proof of our main theorem, before turning to the proof of that proposition. So we consider a ${\rm PO}(1,n)^{\circ}$-invariant decomposition 
$$\HR \otimes \C = U_1\oplus U_2 \oplus W,$$
as in the second case of Proposition~\ref{dichotomy}. We apply Proposition~\ref{totallyreal} to the complex representation of ${\rm PO}(1,n)^{\circ}$ in $U_{1}$. This yields a totally real subspace $V_{1}\subset U_{1}$, whose complexification equals $U_{1}$ and such that the restriction of $B$ to $U_1$ is the Hermitian extension of a strongly nondegenerate bilinear form on $V_1$ of signature $(1,\infty)$. Consider the projections $\HR\otimes \C \to U_1$ with kernel $U_2\oplus W$ and $U_1\to V_1$ with kernel $iV_1$. Let us call $\pi$ the restriction of the composition of these projections to $\HR$, so that $\pi$ is an $\mathbb{R}$-linear map from $\HR$ to $V_1$. 
Note that we can suppose that $\pi\neq 0$, otherwise, we pick the map $Im: U_1\to V_1$ instead of $Re:U_1\to V_1$. Now, $\pi$ is injective since $ker(\pi)$ is a proper closed invariant real subspace of $\mathscr{H}$ and this implies that $ker(\pi)=\{0\}$ by the irreducibility of $\varrho$. By Proposition~\ref{pointwf}, there exists a $2$-dimensional subspace $P$ of $\HR$ wich is pointwise fixed by $K$. The space $\pi(P)$ is contained in $V_1^{K}$ and, by injectivity of $\pi$, ${\rm dim}_{\mathbb{R}}{V_1}^{K}\geq 2$. But this implies that ${\rm dim}_{\mathbb{C}}{U_1}^{K}\geq 2$, a contradiction since the complex representation in $U_1$ is irreducible.

We now turn to the proof of Proposition~\ref{totallyreal}. 

\noindent {\it Proof of Proposition~\ref{totallyreal}.} We call $W$ the complex Hilbert space underlying the representation $\alpha$, and $C$ the Hermitian form of signature $(1,\infty)$ on $W$. Let $\mathbb{H}^{\infty}_{\mathbb{C}}$ be the associated infinite dimensional hyperbolic space, thought of as the space of positive lines in the projective space $\mathbb{P}(W)$. We denote by
$$c : \HHC \times \HHC \times \HHC \to (-\frac{\pi}{2},\frac{\pi}{2})$$
the Cartan argument (see e.g.~\cite{burgeriozzi,monod2018} for the definition). This is an alternating cocycle, which is invariant under the diagonal action of the group ${\rm PU}(1,\infty)$. The representation $\alpha$ defines a continuous isometric action (also denoted $\alpha$) of ${\rm PO}(1,n)^{\circ}$ on $\HHC$. We choose an $\alpha$-equivariant continuous map
$$f : \mathbb{H}^{n}_{\mathbb{R}}\to \HHC.$$
It is classical that any ${\rm PO}(1,n)^{\circ}$-invariant alternating cocycle defined on $\mathbb{H}^{n}_{\mathbb{R}}$ vanishes if $n\ge 3$, see for instance~\cite{monod2018} for a proof. Hence we have $c(f(x),f(y),f(z))=0$ for $x$, $y$, $z$ in $\mathbb{H}^{n}_{\mathbb{R}}$. But this implies that the image of $f$ is contained in a closed totally real subspace of $\HHC$, see for instance Lemma 2.6 in~\cite{monod2018} for a proof and~\cite[p. 469]{burgeriozzi} for the definition of totally real subspaces of $\HHC$. The intersection $N$ of all closed totally real subspaces of $\HHC$ containing the image of $f$ is ${\rm PO}(1,n)^{\circ}$-invariant. The space $N$ is the intersection of $\mathbb{H}^{\infty}_{\C}$ with $\mathbb{P}(V)$ where $V\subset W$ is a totally real subspace, necessarily $\PON^{\circ}$-invariant. Since $\alpha$ is irreducible we must have that $W=V\oplus iV$ and the $\PON^{\circ}$-action on $V$ must be irreducible. Finally, it is easy to see that the restriction of $C$ to $V$ is strongly nondegenerate of signature $(1,\infty)$. This concludes the proof.\hfill $\Box$


\bigskip
\bigskip
\begin{small}
\begin{tabular}{llll}
Pierre Py & & & Arturo S\'anchez\\
IRMA & & & Instituto de Matem\'aticas\\
Universit\'e de Strasbourg \& CNRS & & & Universidad Nacional Aut\'onoma de M\'exico\\
67084 Strasbourg, France & & & Ciudad Universitaria, 04510 M\'exico DF, M\'exico\\
ppy@math.unistra.fr & & & arturo.sanchez@matem.unam.mx\\    
\end{tabular}
\end{small}
 
 \end{document}